\documentclass[final]{siamltex}
\usepackage{amsfonts}
\usepackage{epsfig}

\title{Optimizing Schr\"odinger functionals using Sobolev gradients: Applications 
to Quantum Mechanics and Nonlinear Optics}

\author{
  Juan Jos\'e Garc\'{\i}a-Ripoll
  \thanks{Universidad de Castilla-La Mancha, Departamento de Matem\'aticas,
    E.T.S.I.  Industriales, Avd. Camilo Jos\'e Cela, 3, Ciudad Real, E-13071,
    Spain. ({\tt jjgarcia@ind-cr.uclm.es})}
  \and 
  V\'{\i}ctor M. P\'erez-Garc\'{\i}a
  \thanks{Universidad de Castilla-La Mancha, Departamento de Matem\'aticas,
    E.T.S.I. Industriales, Avd. Camilo Jos\'e Cela, 3, Ciudad Real, E-13071,
    Spain. ({\tt vperez@ind-cr.uclm.es}).}
}

\begin{document}

\maketitle

\begin{abstract}
  In this paper we study the application of the Sobolev gradients technique to
  the problem of minimizing several Schr\"odinger functionals related to timely
  and difficult nonlinear problems in Quantum Mechanics and Nonlinear Optics.
  We show that these gradients act as preconditioners over traditional choices
  of descent directions in minimization methods and show a computationally
  inexpensive way to obtain them using a discrete Fourier basis and a Fast
  Fourier Transform.  We show that the Sobolev preconditioning provides a great
  convergence improvement over traditional techniques for finding solutions
  with minimal energy as well as stationary states and suggest a generalization
  of the method using arbitrary linear operators.
\end{abstract}



\pagestyle{myheadings}
\thispagestyle{plain}
\markboth{J. J. GARC\'IA-RIPOLL AND V. M. P\'EREZ-GARC\'IA}{OPTIMIZING 
SCHR\"ODINGER FUNCTIONALS WITH SOBOLEV GRADIENTS}
\section{Introduction}

The observation of Nature reveals that many unforced continuous systems tend to
accommodate into stationary configurations, in which the distributions of mass,
charge, velocity, etc, do not change throughout time. In the language of
mathematical modeling all configurations are represented by points of a
certain space of functions, $\psi(\mathbf{x}) \in W$, while the tendency of the
system to lie into any of these states is given by a functional, $E(\psi):\, W
\rightarrow\mathbb{R}$, {\em the energy}, whose minima are precisely those
stationary ``states''. For this reason it is possible to see many physical
problems written as variational principles of the type ``find $\psi \in W$ such
that $E(\psi):\, W\rightarrow\mathbb{R}$ achieves a minimum on $W$''. In most
situations the functional to be minimized has a dependence on $\psi$ of the form
\begin{equation}
  \label{functional}
  E(\psi)=\int f\left( \nabla \psi(\mathbf{x}),\psi(\mathbf{x})\right) d^nx.
\end{equation}

However, a complete analytical description of the minima of the functional
$E(\psi)$ is usually not possible. In this paper we will introduce several
techniques for performing this study numerically, focusing on the minimization
of $E(\psi)$ subject to physical constraints.

From a practical point of view, this problem is similar to that of finding the
minima of a real function defined over a finite-dimensional space, such as
$\mathbb{R}^n$.  First, a definition of derivative of the functional, $\nabla
E(\psi)$, must be chosen. If the domain of the functional $W$ is equipped with
some scalar product we may use the Frechet derivative which is given by a first
order expansion of the functional around a function $\psi$
\begin{equation}
  \label{gradient}
  E(\psi+\delta)=E(\psi)+\langle\delta,\nabla E(\psi)\rangle
  +\langle\nabla E(\psi),\delta\rangle
  +O\left( \left\Vert\delta\right\Vert^2\right) .
\end{equation}

The \emph{critical points,} $\psi_c$, are defined as the points where the first
order variation of the functional vanishes for any perturbation $\delta$. That
is, the derivative vanishes in a weak sense
\begin{equation}
  \label{critical-point}
  \langle \delta,\nabla E(\psi_c) \rangle =0,\quad \forall \delta.
\end{equation}

Just like in the finite-dimensional case it is possible to show that any
minimum of the functional must be also a critical point. Thus a common approach
is to solve Eq. (\ref{critical-point}) and verify {\em a posteriori} which
solutions are actually minima of the functional. If $W = L^2(\mathbb{R}^n)$
this procedure gives us the well known Euler-Lagrange equations of the problem,
which is a partial derivatives equation (PDE)
\begin{equation}
  \label{Lagrange}
  \frac{\partial f}{\partial \bar{\psi}}-\nabla \cdot \frac{\partial f}{\partial \nabla \bar{\psi}}=0.
\end{equation}

However, we have no guarantee to reduce the complexity of the problem, as it is
by no means trivial to solve Eq. (\ref{Lagrange}). Also we are likely to obtain
more solutions than we actually need, since not only minima, but also maxima
and saddle node points will satisfy the Lagrange equations.

To avoid these problems some other methods are used which aim at finding the
minima of the functional directly, constructing minimizing sequences,
$\{\psi_i\}$, whose limit is a minimum of the functional:
$\psi_c=\lim_{i\rightarrow \infty }\psi_i$.  These methods will be discussed in
the following sections.

The outline of this paper is as follows. In \S \ref{sect-theory} we recall the
definition of Sobolev gradients as given in \cite{Neu98,Neu97}. We derive a
formal solution to the problem of finding these gradients which is based on the
inversion of a positive hermitian operator. In \S \ref{sect-Fourier} we derive
an explicit expression for the Sobolev gradients in the trigonometric Fourier
basis and comment on its implementation using Fast Fourier Transforms (FFT).

In \S \ref{sect-bose} and \S \ref{sect-optics} we apply the previous tools to
two physical problems. Using descent techniques with Sobolev gradients over
Fourier spaces we will find the ground states of a Bose-Einstein condensate in
a rotating magnetic trap and the excited states for coupled laser beams
propagating through a nonlinear medium. Both physical systems are modeled by
nonlinear equations of Schr\"odinger type and present difficulties when
traditional minimization techniques are used. We comment on the great
improvements that are achieved using Sobolev gradients.  Finally in \S
\ref{sect-final} we summarize our results and offer some conclusions.

\section{Sobolev gradients}
\label{sect-theory}

\subsection{Direct solutions of the variational problem}

There are two traditional approaches to the problem of finding the minima of a
functional using a discrete basis. The first one expands the unknown solution
using a Fourier basis $\{\phi_k\}$, $\psi=\sum_kc_k\phi_k,$ and defines a new
functional over the finite--dimensional space
\begin{equation}
  \label{finite-dimensional}
  E(\{c_k\})\equiv E\left(\sum c_k\phi_k\right).
\end{equation}
The functional is then minimized using methods which are well known from the
the domain of finite--dimensional problems, e.g. Newton's method or nonlinear
conjugate gradient.

This procedure is quite straightforward and there is a huge amount of
literature and tools which can be immediately applied to
(\ref{finite-dimensional}).  However, for some types of problems one has to
deal with highly nonlinear algebraic equations with many of terms on each
equation, something which is computationally too expensive to work with.

The second approach involves what is known as \emph{descent techniques}.  The
idea is to manipulate the original functional (\ref{functional}) building an
analytic equation for a minimizing trajectory in the target space $W$. This
equation is then discretized and solved on a suitable basis. In the
\emph{continuous steepest descent} version, the trajectory $\psi(t):\,
\mathbb{R} \rightarrow W$ is continuous and defined by a PDE which involves the
gradient of the functional (\ref{gradient})
\begin{equation}
  \label{cont-sd}
  \frac{\partial \psi}{\partial t}=-\nabla E.
\end{equation}
The \emph{discrete steepest descent} technique is computationally cheaper since
instead of requiring an integrator for the PDE it constructs a discrete
succession of estimates to the minimum, $\{\psi_{k+1}=\psi_k+\lambda_k\nabla
E(\psi_k)\}$, by locally minimizing $E(\psi_k+\lambda\nabla E)$ with respect
to the real parameter $\lambda$.

In this paper we will only deal with descent techniques. The first and most
important reason is that we will work with the definition of $\nabla E(\psi_k)$
trying to improve its convergence. As it was already shown in \cite{Neu97},
this work pays off: a good choice of the gradient improves convergence by
several orders of magnitude. The second motivation is that by focusing on the
gradient our algorithms will be essentially independent on the descent method,
which leaves space for further improvement. For instance one might apply these
techniques to a nonlinear conjugate gradient method ---which dynamically
adjusts the search direction, $d_k\equiv \nabla E(\psi_k)$, with an estimate
that takes into account the history of the evolution---. Finally, by working
with Eq.  (\ref{cont-sd}) or its discrete version we avoid the complex
nonlinearities that arise in other methods, and we will have a more ample
choice of Fourier basis to work with.

\subsection{Ordinary gradients}

It is customary in the literature to work in spaces which are equipped with an
$L^2$ scalar product and its corresponding norm
\begin{eqnarray}
  \langle \psi,\phi \rangle _{L^2} & \equiv  & \int \bar{\psi}\phi ,\label{L2-product} \\
  \left\Vert \psi\right\Vert^2_{L^2} & \equiv  & \int |\psi|^2.\label{L2-norm}
\end{eqnarray}
If one does so and works with Eq. (\ref{gradient}) then the formal definition
of the gradient is Lagrange's one
\begin{equation}
  \label{ordinary-gradient}
  \nabla E(\psi)=\frac{\partial E}{\partial \bar{\psi}}-\nabla \frac{\partial E}{\partial \left( \nabla \bar{\psi}\right) }.
\end{equation}
We will refer to this definition as the ``ordinary'' gradient, to distinguish
it from the different definitions that we will derive below.

\subsection{Sobolev gradients}

Following the ideas from \cite{Neu97} we will move our problem to a different
space which is the Sobolev space of functions, such that $\psi$ and its
derivatives, $\nabla \psi$, have a well defined $L^2$-norm:
\begin{equation}
  \mathbb{H}^1\equiv \{\psi/\psi,\nabla \psi\in L^2\}.
\end{equation}
This Sobolev space will also be equiped with a scalar product and a norm
\begin{eqnarray}
  \langle \psi,\phi \rangle  & \equiv  &
  \int \left[ \bar{\psi}(\mathbf{x})\phi(\mathbf{x})+\nabla \bar{\psi}(\mathbf{x})\cdot \nabla \phi(\mathbf{x})\right] d^nx,\\
  \left\Vert \psi\right\Vert^2 & \equiv  &
  \int \left[ |\psi|^2+|\nabla \psi|^2\right] d^{n}x.
\end{eqnarray}

To obtain a new explicit expression for the gradient of the functional in the
Sobolev space we will follow a less rigorous derivation than in \cite{Neu98}.
Performing a first order expansion of our functional around a trial state
$\psi$ we obtain
\begin{equation}
E(\psi+\varepsilon \delta )=E(\psi)+\varepsilon \int \overline{\left( \delta
    ,\nabla \delta \right) }\left( \begin{array}{c}
    \frac{\partial E}{\partial \bar{\psi}}\\
    \frac{\partial E}{\partial \nabla \psi}
\end{array}\right) +c.c.+O(\varepsilon^2).
\end{equation}

We have to turn this expression into something like Eq. (\ref{gradient}).  This
means that we have to find some $\phi$ such that
\begin{equation}
  \int \left[ \bar{\delta }\frac{\partial E}{\partial \bar{\psi}}+\nabla
    \bar{\delta }\frac{\partial E}{\partial (\nabla \bar{\psi})}\right]
  +c.c=\int \left[ \bar{\delta }\phi +\nabla \bar{\delta }\nabla \phi \right]
  +c.c.
\end{equation}
If we integrate by parts and impose that this equality be satisfied for all
perturbations, $\delta$, the problem has a formal solution which is given by a
Lagrange equation
\begin{equation}
  \left( 1-\triangle \right) \phi =\frac{\partial E}{\partial \bar{\psi}}-\nabla\frac{\partial E}{\partial (\nabla \bar{\psi})}.
\end{equation}
In consequence, our formal expression for the Sobolev gradient of $E(\psi)$
finally reads
\begin{equation}
  \label{sobolev-gradient}
  \nabla_SE\equiv \left( 1-\triangle \right) ^{-1}\nabla E.
\end{equation}
Here $\nabla_SE$ stands for the Sobolev gradient, $\nabla E$ is the ordinary
one, and $(1-\triangle)^{-1}$ represents the inverse of a linear and strictly
positive definite operator.

\section{Sobolev gradients on Discrete Fourier spaces}
\label{sect-Fourier}

In the rest of the paper we will work with functions which are defined over a
$d$-dimensional rectangular volume $\Omega \equiv \{x\in \Pi_i[a_i,b_i]\}$ with
side lengths given by $L_i=b_i-a_i$.  We customarily define an orthogonal set
of basis functions, $\phi_n=e^{ik_nx}$, over $\Omega$ where $ \left\{k_n=2\pi
  \left(\frac{n_1}{L_1},\ldots,\frac{n_d}{L_d}\right),\, n_i\in Z\right\}$

It is well-known that it is possible to expand any continuous function $f(x)$
with periodic boundary conditions using this basis
\begin{equation}
  f(x)=\sum^{+\infty}_{n=-\infty}\hat{f}_n\phi_n(x)
\end{equation}
where
\begin{equation}
  \label{coeff-Fourier}
  \hat{f}_n=\frac{1}{V}\int_{\Omega}\bar{\phi}_n(x)f(x).
\end{equation}
In the previous formula $V=\Pi_iL_i$ is the volume of $\Omega$ and arises
because of the lack of normalization of the basis functions, a common practice
which saves some computation time.

To discretized the problem we will work within the set of functions sampled over
a set of evenly spaced points from $\Omega$, $\{x_n=(n_1h_1,\ldots,n_dh_d),\,
n_i=0, \ldots ,N_i-1\}$. Here $h_i$ represents the spacing along the $i$-th
dimension, $n$ is a vector of non--negative integers and we denote a sampled
function by an index, as in $f_n\equiv f(x_n)$.

Due to this discretization our previous Fourier basis is now redundant. It can
be reduced to a finite subset of functions which may represent any sampled
function. These functions are given by $\left\{k_n=2\pi \left(
    \frac{n_1}{L_1},\ldots ,\frac{n_d}{L_d}\right),\,
  n_i=-M_i+1,\ldots,M_i\right\},$ where $M_i=[N_i/2]$ is the integer quotient
of $N_i$ divided by two. In this basis a sampled function is given by an
expansion which reads
\begin{equation}
  \label{discrete-fourier}
  f_m=\sum_n\hat{f}_n\phi_n(x_m),
\end{equation}
and which makes use of the same coefficients as Eq. (\ref{coeff-Fourier}).

The advantage of the finite Fourier basis over other approaches is that it
provides an approximant of any function whose error is of order
$O(L_i/N_i)^{p}$ where $p$ is the maximum differentiability of the sampled
function. Furthermore, there is a numerically efficient method known as Fast
Fourier Transform which allows one to compute the Fourier coefficients up from
the sampled function, $f_n\rightarrow \hat{f}_m$, and vice versa.

Solving Eq. (\ref{sobolev-gradient}) numerically in a discrete Fourier basis is
simple. Let us say that we have computed the ordinary gradient and that its
sampled version has some Fourier coefficients
\begin{equation}
  \nabla E(\mathbf{x}_n)=\sum \hat{e}_m\phi_m(x_n).
\end{equation}
and let us assume that there exists a certain solution to Eq.
(\ref{sobolev-gradient}) and that it has another discrete Fourier expansion
\begin{equation}
  \nabla_sE(\mathbf{x}_n)=\sum \hat{s}_m\phi_m(x_n).
\end{equation}
Then by virtue of Eq. (\ref{sobolev-gradient})
\begin{equation}
  \hat{s}_m=\frac{\hat{e}_m}{1+k^2_m}.
\end{equation}
That is, in the sampled space the Sobolev gradient represents a preconditioning
of the ordinary gradient, such that the most oscillating modes are more
attenuated.  Furthermore, due to this very simple expression computing the
Sobolev preconditioning is computationally cheap and involves only minor
changes to existing computer codes based on Fourier transforms.

\section{Applications to Quantum Mechanics}
\label{sect-bose}

\subsection{The problem}

In this section we will apply the Sobolev gradients technique to a timely
problem from Quantum Physics. The system that we will study is a dilute gas of
bosonic atoms which are cooled down to ultralow temperatures at which their
dynamics become synchronized. When the temperature is low enough, the gas or
``condensate'' may be described using a single complex wave function,
$\psi(\mathbf{x},t)$ which is ruled by the so called Gross-Pitaevskii equation,
a type of Nonlinear Schr\"odinger equation that for the Bose gas in a rotating
trap with adimensional units reads
\begin{equation}
  \label{GPE-rot}
  i\partial_t\psi(\mathbf{x},t)=
  \left[ -\frac{1}{2}\triangle +V(\mathbf{x})+g|\psi(\mathbf{x},t)|^2-\Omega L_z\right]\psi(\mathbf{x},t).
\end{equation}
Here $g\in \mathbb{R}^+, \Omega \in \mathbb{R}$, $L_z = i\left(
  x_1\partial_2-x_2\partial_1\right)$ is an hermitian operator whose expected
value $\langle L_z\rangle =\int \bar{\psi}L_z\psi,$ represents the
\emph{angular momentum} of the condensate along an axis of the trap.

There is conserved quantity associated to equation (\ref{GPE-rot}) which is
called the \emph{energy functional} of the condensate
\begin{equation}
  \label{GPE-energy}
  E(\psi)=
  \frac{1}{2}\int \left\{|\nabla \psi|^2+\bar{\psi}\left[ V(\mathbf{x})+\frac{1}{2}g|\psi|^2-\Omega L_z\right] \psi\right\} d^nx.
\end{equation}
Our objective in this part of the work will be to \emph{find the solutions
  which are the minima of the energy subject to a restriction of the
  $L^2$-norm}
\begin{equation}
  \label{norm-restriction}
  \int |\psi|^2\equiv N.
\end{equation}
The particular value of $N$ is imposed by the experimental conditions.
Furthermore, without this restriction the absolute minimum of the energy is
always reached at the trivial solution $\psi=0$.

In Quantum Mechanics the variational formulation of the problem is
traditionally converted into a Lagrange equation (\ref{Lagrange}), which is
nothing but the Gross-Pitaveskii equation for so called {\em stationary
  states}. In short the word {\em ``stationary''} refers to solutions of the
type
\begin{equation}
  \psi(\mathbf{x},t)=e^{-i\mu t}\phi_{\mu} (\mathbf{x}).
\end{equation}
The pair $\{\mu, \phi_{\mu} (\mathbf{x})\}$ satisfies a nonlinear eigenvalue
problem
\begin{equation}
  \label{GPE-stat}
  \mu \phi_{\mu}(\mathbf{x})=
  \left[ -\frac{1}{2}\triangle +V(\mathbf{x})+g|\phi_{\mu} (\mathbf{x})|^2-\Omega L_z\right] \phi_{\mu} (\mathbf{x}).
\end{equation}
Due to the difficulty of solving directly the problem (\ref{GPE-stat}) we
will try to find a direct solution to the variational problem.

\subsection{Numerical methods: imaginary time evolution}

We can search the minima of (\ref{GPE-energy}) using descent techniques
modified to account for the restriction on the norm (\ref{norm-restriction}).
The first way to do this is to use a version of the continuous steepest descent
which is known as \emph{imaginary time evolution} \cite{imaginary-time}. We can
summarize this method with the following set of equations
\begin{eqnarray}
  \nu (\mathbf{x},\tau) & = & \sqrt{\frac{N}{\left\Vert \sigma \right\Vert _{L^2}^2}}\sigma (\mathbf{x},\tau),\label{bec-imaginary-time-1}\\
  \frac{\partial \sigma }{\partial \tau}(\mathbf{x},\tau) & = & -\nabla E(\nu ),\label{bec-imaginary-time} \\
  \nabla E(\nu ) & = & \left[ -\frac{1}{2}\triangle +V(\mathbf{x})+g|\nu |^2-\Omega L_z\right] \nu .
\end{eqnarray}
Here we see that $\nu(\mathbf{x},\tau)$ evolves continuously maintaining a
fixed $L^2$--norm $N$ and following the direction of decreasing energy given by
$\nabla E.$ Indeed it is easy to show that $\frac{\partial }{\partial
  \tau}\left[ E(\nu (\mathbf{x},\tau))\right] \leq 0$.  Hence, the limit given
by
\begin{equation}
  \label{method-convergence}
  \phi(\mathbf{x})=\lim_{\tau \rightarrow \infty }\nu(\mathbf{x},\tau)
\end{equation}
is at least a critical point of the energy, if not a minimum\footnote{%
  As it is common with these local minimization procedures, it remains the
  problem that iterations may be trapped in a critical point which is not a
  minimum.  Linear stability analysis may then be applied to check the validity
  of the solution.}

From a practical point of view, the simplest way to find the minimizer using
imaginary time evolution is to repeatedly integrate Eq.
(\ref{bec-imaginary-time}) for a very short time, $\Delta t$, apply Eq.
(\ref{bec-imaginary-time-1}) for the newly found $\sigma(\mathbf{x},t+\Delta
t)$ and use the new estimate for $\nu$ to redefine $\sigma\equiv\nu$ and repeat
the procedure until convergence.  This way one avoids the problem that
according to Eq.(\ref{bec-imaginary-time}) the norm of $\sigma$ may grow
indefinitely. The same consideration applies to the remaining methods that we
will present here.

Although the method from Eq. (\ref{bec-imaginary-time}) was derived using an
ordinary gradient, nothing prevents us from applying our Sobolev
preconditioning and all results should be still valid. If we do so, our new
equations are
\begin{eqnarray}
  \nu (\mathbf{x},\tau) & = & \sqrt{\frac{N}{\left\Vert \sigma \right\Vert _{L^2}^2}}\sigma (\mathbf{x},\tau),\\
  \frac{\partial \sigma }{\partial \tau}(\mathbf{x},\tau) & = & -\nabla_SE(\nu ),\label{bec-imaginary-time-sobolev} \\
  \nabla_SE(\nu ) & = & \left( 1-\triangle\right)^{-1}\left[ -\frac{1}{2}\triangle +V(\mathbf{x})+g|\nu|^2-\Omega L_z\right] \nu,
\end{eqnarray}
and the critical point is still given by Eq. (\ref{method-convergence}).

\subsection{Numerical methods: minimization of the free energy}

While the imaginary time evolution is easy to understand and to implement, the
fact that it performs the descent over a path of functions with a certain norm
makes it too restricted and sometimes too slow. A different approach is to
define a new functional called the \emph{free energy} with a Lagrange
multiplier that takes care of the restriction on the norm. In Quantum Mechanics
this free energy is usually defined as
\begin{equation}
 \label{old-free-energy}
  F_{QM}(\psi)=E(\psi)-\mu N(\psi),
\end{equation}
because it preserves the linear form of the equations. Since our equations are
already nonlinear we define a free energy functional more conveniently as
\begin{equation}
\label{new-free-energy}
  F(\psi)=E(\psi)+\frac{1}{2}\left(N(\psi)-\lambda \right)^2.
\end{equation}

First and most important it is not difficult to show that any absolute or
relative minimum of $F(\psi)$ is also a minimum of $E(\psi)$ subject to Eq.
(\ref{norm-restriction}) with a nonlinear eigenvalue given by $\mu
=N(\psi)-\lambda $.

Secondly, as we will prove in Appendix \ref{free-energy-properties}, $F(\psi)$
must have at least one finite norm minimum, something which cannot be easily
assured for $F_{QM}$.

The practical advantage of our new functional consists in that fixing $\lambda$
and $\Omega$ we can perform a continuous descent over the whole domain of
$F(\psi)$ without renormalizing the solution on each iteration ---i.e., the
search space is larger.  The new equation that we must integrate is thus
\begin{equation}
  \frac{\partial \nu }{\partial \tau}(\mathbf{x},\tau)  =  -\nabla F(\nu ) = - \left[ -\frac{1}{2}\triangle +V(\mathbf{x})+g|\nu |^2-\Omega L_z\right] \nu .\label{bec-free-energy}
\end{equation}
Using a Sobolev preconditioning Eq. (\ref{bec-free-energy}) becomes
\begin{equation}
  \frac{\partial \nu }{\partial \tau}(\mathbf{x},\tau) = -\nabla_SF(\nu )=-\left( 1-\triangle\right)^{-1}\left[ -\frac{1}{2}\triangle +V(\mathbf{x})+g|\nu |^2-\Omega L_z\right] \nu.\label{bec-free-energy-sobolev} 
\end{equation}
In both cases the actual solution is still given by the limit of Eq.
(\ref{method-convergence}).

\subsection{Numerical results}

Up to this point we have shown four different numerical methods, two of them
incorporating a Sobolev preconditioning [Eqs.
(\ref{bec-imaginary-time-sobolev}) and (\ref{bec-free-energy-sobolev})] and two
without it [Eqs.  (\ref{bec-imaginary-time}) and (\ref{bec-free-energy})].  We
have compared the efficiencies of these methods for several test situations.
The details of our study are as follows.

All methods have been implemented using the discrete Fourier basis mentioned
above. This applies both to the calculation of derivatives and to the
application of the Sobolev preconditioning. We restricted our simulations to
two-dimensional problems with radially symmetric potential,
$V(\mathbf{x})=\frac{1}{2}|\mathbf{x}|^2$, over a grid of 128$\times$128
points. Practical experience with more complex problems shows that our results
generalize to higher dimensionality and denser grids.

Due to the requirements of the method, both variants of imaginary time
evolution [Eq. (\ref{bec-imaginary-time}) and
(\ref{bec-imaginary-time-sobolev})] are integrated using a Runge-Kutta-Fehlberg,
 where the tolerance is adapted as the solution
converges to its target. On the other hand, instead of performing a continuous
descent for both variants of the free energy descent [Eq.
(\ref{bec-free-energy}) and (\ref{bec-free-energy-sobolev})], we found it more
convenient and faster to perform a discrete steepest descent.

All programs have been implemented using the tensor-algebra environment called
Yorick \cite{Yorick}, an interpreted environment which is capable of fast
numerical computations and which can be equipped with the FFTW library
\cite{FFTW}. Execution times are given for a Digital Personal workstation
500au, but the programs run equally well on modest personal computers with
Pentium-II processors and less than 64 Mb of memory.

As test case we have considered three situations. Case A is the simplest one
and corresponds to a stationary trap $\Omega =0$, an intense nonlinearity
$g=100$ and starting the minimization with a radially symmetric Gaussian of
unit width, that is, $\psi_0 \propto e^{-|\mathbf{x}|^2}$, a shape which is
similar to that of the true ground state.

Cases B and C involve a rotating condensate with $\Omega=0.6$ and $g=100$.
Under these conditions the functional achieves both an absolute minimum, which
is the same that found in case A, and a local minimum. The local minimum is a
solution of vortex type, a topological defect, whose behavior near zero is
$\psi \propto (x_1+ix_2)/|\mathbf{x}|^2$.

For this reason we designed two test cases with different initial conditions.
Case B starts with a Gaussian profile with a centered vortex or $\psi_0 \propto
|\mathbf{x}|e^{-|\mathbf{x}|^2}(x_1+ix_2)/|\mathbf{x}|^2.$ In this case all
methods are trapped on the local minimum with the centered vortex.

Finally, a third set of simulations, case C, uses the same parameters
$\{\Omega=0.6, g=100\}$ but starts from a nonsymmetric initial state, $\psi_0
\propto |\mathbf{x}| e^{-|\mathbf{x}|^2}
((x_1-y_1)+i(x_2-y_1))/|\mathbf{x-y}|^2$ which is close to the vortex solution
but belongs to the basin of attraction of the ground state.  Here all
minimization methods require more computational work since they must find a
path out of the local minimum, and it is precisely in this case where the
differences between methods are best shown.

\begin{figure}
\begin{center}
\epsfig{file=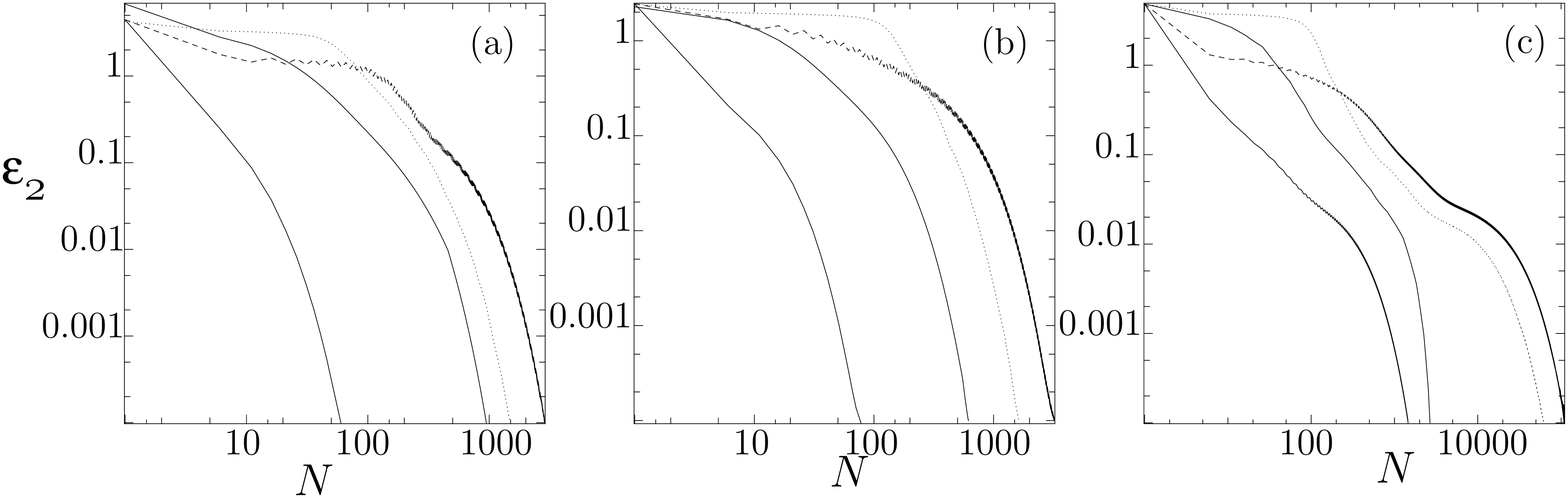,width=13cm}
\end{center}
\caption{Evolution of error, $\varepsilon_2\equiv\Vert\psi-\psi_{exact}\Vert_2$,
  through different minimization processes for continuous steepest descent with
  Sobolev preconditioning (lower solid line) and without it (dashed line), and
  for imaginary time evolution with Sobolev preconditioning (upper solid line)
  and without it (dotted line).  Plots (a) to (c) correspond respectively to
  the cases A,B, and C described in the text. Both axes, error and number of
  iterations, are in a logarithmic scale.}
\end{figure}

In Table \ref{table-bec} we summarize the results of the simulations. In case A
it is apparent that the Sobolev preconditioning has a positive influence over
convergence, with an astonishing result of 55 steps for the steepest descent
with free energy. A similar behavior is found in case B.

\begin{table}
  \vspace{0.3cm}
  \begin{center}
    \begin{tabular}{|c|c|c|c|c|c|}
      \hline
      \multicolumn{2}{|l|}{Methods} & IT &  ITS   & FE  & FES \\ \hline\hline
      Case A & Iterations           & 1320  & 945  & 2850   & 55  \\ \cline{2-6}    
             & Time (s)             &  416  & 371  & 285    & 13 \\ \hline
      Case B & Iterations           & 1630  & 615  & 3210   & 320 \\ \cline{2-6}
             & Time (s)   &  468  & 242  & 75     & 9 \\ \hline  
      Case C & Iterations & 64195 & 2665 & 108505 & 1455 \\ \cline{2-6}
             & Time (s)   & 19863 & 1165 & 10861  & 168 \\ \hline
    \end{tabular}
  \end{center}
  \vspace{0.3cm}
  \caption{\label{table-bec}
    Iterations and computation time for each minimization method: imaginary time
    without (IT) and with (ITS) Sobolev preconditioning and free energy without
    (FE) and with (FES) Sobolev preconditioning. Shown are results for the initial
    data described in the text (cases A, B and C).}
\end{table}

In case C, the Sobolev preconditioning enhances convergence speed by two orders
of magnitude. An intuitive explanation of why the steepest descent with a
Sobolev gradient takes less steps to converge will be discussed in detail in
Appendix \ref{generalized-gradients}.

\section{Applications to Nonlinear Optics\label{sect-optics}}

\subsection{The model}

In this section we will consider a model for a pair of incoherently interacting
light beams. To be precise we will study the light field of each beam,
$u(\mathbf{x},t)$ and $w(\mathbf{x},t)$, propagating through a weakly nonlinear
saturable optical medium.  This system may be modeled by the Cauchy problem
\begin{eqnarray}
i\partial_tu & = & -\triangle u+\frac{u}{1+\kappa (|u|^2+|w|^2)},\label{beams-equations1} \\
i\partial_tw & = & -\triangle w+\frac{w}{1+\kappa (|u|^2+|w|^2)}.\label{beams-equations2}
\end{eqnarray}
for the complex functions $u,w:\mathbb{R}^2\times \mathbb{R}^+ \rightarrow
\mathbb{C}$, which vanish at infinity and satisfy the initial data
$u(\mathbf{x},0)=u_0(\mathbf{x})$ and $w(\mathbf{x},0) = w_0(\mathbf{x})$.
Here $\kappa\in \mathbb{R}^+$, $-\triangle = \frac{\partial^2}{\partial x^2} +
\frac{\partial^2}{\partial y^2}$ is the Laplacian operator which accounts for
the diffraction of light and the nonlinear term $(1+|u|^2+|w|^2)^{-1}$ models
the saturable interaction among the beams.

Let us define the two component vector
\begin{equation}
\tilde{U}(\mathbf{x},t)=\left( \begin{array}{c}
    u(\mathbf{x},t)\\
    w(\mathbf{x},t)
\end{array}\right) ,
\end{equation}
then, the energy functional for the system reads
\begin{equation}
E(\tilde{U})=\int \left[ -\tilde{U}^{\dagger }\triangle \tilde{U}+G(|\tilde{U}|^2)\right] d^nx.
\end{equation}
where $G(\rho)= \frac{1}{\kappa^2}\left(\ln \left( 1+\kappa\rho\right)
  -\kappa\rho\right)$. The analysis of this section may be generalized to more
general nonlinearities with only minor changes to $G(\rho)$.

\subsection{Stationary solutions}

We are interested on stationary solutions, which are of the form
\begin{equation}
\tilde{U}(\mathbf{x},t)=\left( \begin{array}{cc}
    e^{i\mu_u t} & 0\\
    0 & e^{i\mu_w t}
\end{array}\right) U(\mathbf{x})=e^{iMt}U(\mathbf{x}).
\end{equation}
The equations for the stationary solutions now are of elliptic type
\begin{equation}
\label{beams-stationary-eq}
MU=-\triangle U+G'(|U|^2)U.
\end{equation}
with zero Dirichlet boundary conditions at infinity. Again this formulation
poses a nonlinear eigenvalue problem for the pair $\{ M, U\}$.

The stationary solutions are critical points of the energy functional subject
to a constraint on the $L^2$-norm of each component. That is, defining
\begin{eqnarray}
\label{beams-constraint1}
N_u & = & \int |u|^2d^{n}x,\\
N_w & = & \int |w|^2d^{n}x,\label{beams-constraint2}
\end{eqnarray}
the first order variation of the energy around $U_0$ for fixed
$N_u$ and $N_w$ must be zero
\begin{equation}
\left. \frac{\delta E}{\delta U}\right| _{N_u,N_w}=0.
\end{equation}
The particular fixed values of $\{N_u,N_w\}$ represent the total intensity of
each beam of light.

\subsection{Ground states}

In principle there are many different stationary solutions of Eq.
(\ref{beams-stationary-eq}).  However, if we focus on the ground states
(stationary solutions of minimal energy) we can apply some of the methods that
we mentioned in \S \ref{sect-bose} with minor modifications to account for the
higher dimensionality of the problem.  For instance, one may define a free
energy functional
\begin{equation}
  F(U)=E(U)+\frac{1}{2}\left( N_u-\lambda_u\right)^2+\frac{1}{2}\left(N_w-\lambda_w\right)^2,
\end{equation}
and minimize it using a discrete steepest descent with Sobolev preconditioning.

\begin{figure}
\begin{center}
\epsfig{file=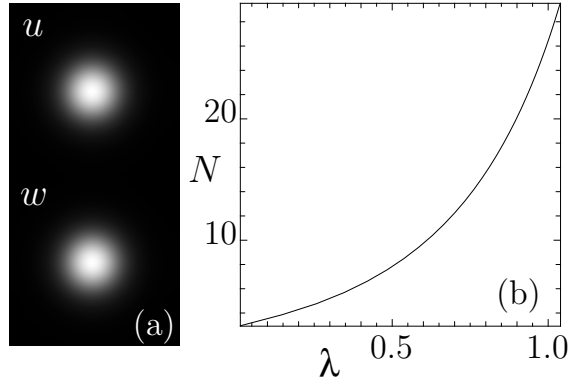,width=7.5cm}
\end{center}
\caption{\label{fig-soliton}
  (a) Density profile of a ground state for $N_u=N_w=30$. (b) Norms 
$N_u=N_w=N$, as a function of the nonlinear Lagrange multipliers
  $\lambda_u=\lambda_w=\lambda$.}
\end{figure}

By performing this minimization using different parameters
$\{\lambda_u,\lambda_w\}$ one obtains nodeless localized solutions for $u$ and
$w$ as shown in Fig.  \ref{fig-soliton}. The precise values of $\lambda_u$ and
$\lambda_w$ determine the norm of each component.

Let us remark that, up to a global factor, the ground state has the same shape in
both envelopes. That is,
\begin{eqnarray}
u_0(\mathbf{x}) & = & N \sqrt{\xi} \, \rho(|\mathbf{x}|,N), \\
w_0(\mathbf{x}) & = & N \sqrt{1-\xi} \, \rho(|\mathbf{x}|,N),\quad \forall\xi\in[0,1]. \nonumber
\end{eqnarray}
The common shape $\rho$ depends on the total intensity $N=N_u+N_w$, and it is
the one that fixes the values of $\mu_u$ and $\mu_w$. For this reason if we had
chosen the traditional definition of the free energy $F_{QM}(U) = E(U) + \mu_u
N_u + \mu_w N_w$ we would have found an infinite degeneracy which prevents
convergence to the desired values of $\{N_u,N_w\}$. This problem is removed by
the use of our nonlinear Lagrange multipliers and the $\{\lambda_u,\lambda_w\}$
parameters.

\subsection{Excited states}
\label{excited-search}

In the field of guided light waves there has been a great interest on the
properties of solutions of Eq. (\ref{beams-stationary-eq}) which are not ground
states, the so called {\em excited states}.  Minimization methods based on the
energy functional [See \S \ref{sect-bose}] cannot be applied to this task since
excited states are not necessarily local minima of the energy. Instead, a
variational principle somehow equivalent to Eq.  (\ref{beams-stationary-eq})
must be defined.

Let us rewrite Eq.  (\ref{beams-stationary-eq}) as the application of a
nonlinear operator
\begin{equation}
  \label{beams-nonlineq}
  f(U)\equiv \left[ -\triangle -M+I(|U|^2)\right] U=0,
\end{equation}
define the error functional 
\begin{equation}
 \label{beams-error}
  F(U)\equiv \int f(U)^{\dagger}f(U)\geq 0
\end{equation}
and take our variational principle to be \emph{``find $U_0$ such that $ F$
  reaches an absolute minimum $F(U_0)=0$.''}

With this principle and a Sobolev preconditioning our descent technique becomes
\begin{equation}
  \label{sobolev-nonlineq}
  \frac{\partial \nu (\mathbf{x},\tau)}{\partial \tau}  = -\nabla_SF(\nu ) = - \left( 1-\triangle \right) ^{-1}\nabla F,
\end{equation}
and it may be proven that $U_0(\mathbf{x}) = \lim _{\tau \rightarrow \infty
  }\nu (\mathbf{x},\tau)$. The ordinary gradient $\nabla F$ reads
\begin{eqnarray}
  G & \equiv  & (-\triangle +I(|U|^2))U,\\
  \nabla F & = & (-\triangle +I(|U|^2))G+I'(|U|^2)\left( U^{\dagger }G+G^{\dagger }U\right) U,
\end{eqnarray}
Its calculation using a discrete Fourier basis is straightforward.

There are several advantages of this approach. The first one is that $F(U)$
makes no distinction between ground and excited states: any stationary state
with the right eigenvalues $\mu_i=M_{ii}$ is a local minimum of this new
functional.  The second one is that we do not need to add any Lagrange
multipliers to $F(U)$ since they are already present in the $M$ operators.
Finally we expect that the number of minima of $F(U)$ be discrete and separated
so as to avoid problems with descent methods being trapped on critical points
that are not minima.  Indeed, it is very easy to know when this accidental
trapping happens since any absolute minimum of $F$ must also be a zero of it
$F(U_0)=0$.

A remarkable feature of the method that we have outlined above is that to
distinguish among the different excited states we have to change both the
ad-hoc eigenvalues, $\{\mu_u,\mu_w\}$, and the initial data of the minimization
method.

\begin{figure}
\begin{center}
\epsfig{file=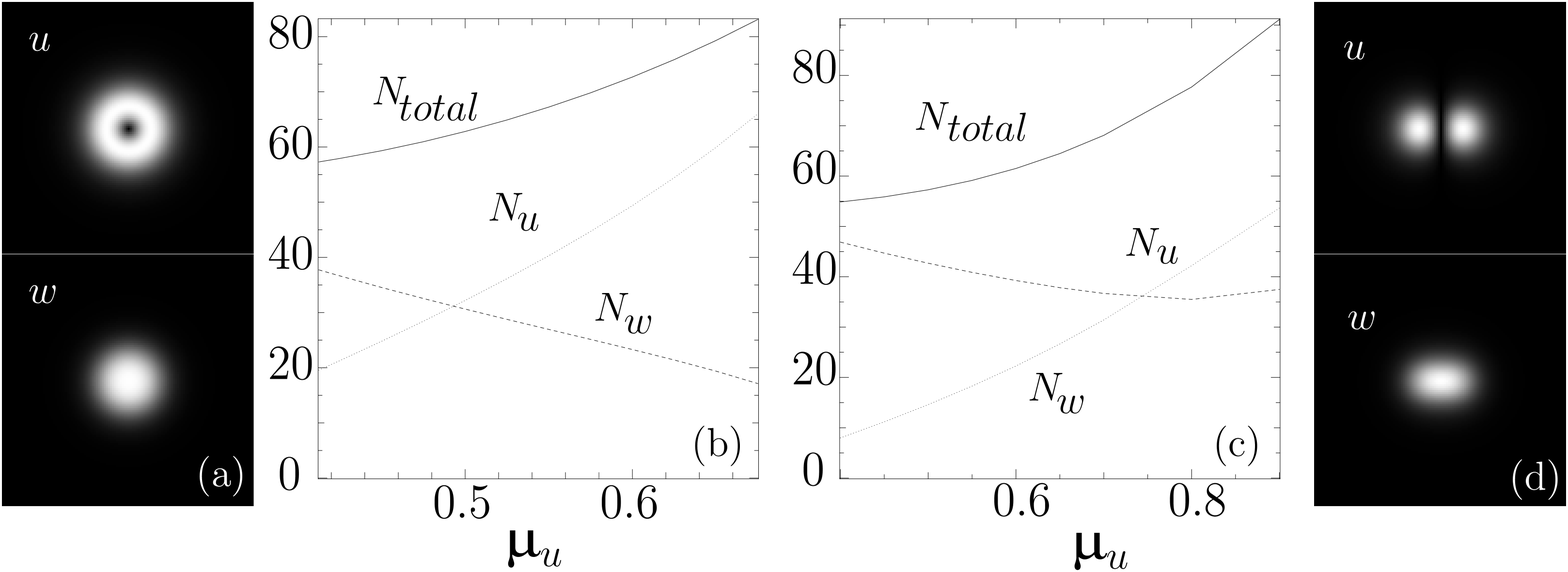,width=13cm}
\end{center}
\caption{\label{fig-vortex}
  (a) Vortex-mode vector solitons and (b) their norms $N_u$,
  $N_w$, $N_{total}=N_u+N_w$, as a function of $\mu_u$ for
  $\mu_w=1$. (d) Dipole-mode vector solitons and (c) their
  norms $N_u$, $N_w$, $N_{total}=N_u+N_w$, as a function of
  $\mu_u$ for $\mu_w=1$. }
\end{figure}

We have concentrated on two types of {\em excited states} of particular
interest for applications: the first one is called \emph{vortex vector soliton}
\cite{Segev}, its features being summarized in Fig.  \ref{fig-vortex}(a-b).
Choosing different initial data for the minimization process we obtain an
second type of excited states which are called \emph{dipole-mode vector
  solitons} \cite{Kivshar}.  An example of these asymmetric solutions is
depicted in Fig.  \ref{fig-vortex}(d). Depending on the parameters $\mu_{u,w}$
we find different norms of the solution. The fact that our method allows the
finding of these nonsymmetric stationary states is very interesting since
conventional approaches to the problem have severe difficulties. In fact, the
application of the Sobolev preconditioning not only greatly enhances the
convergence but it is necessary to get convergence.

\begin{figure}
\begin{center}
\epsfig{file=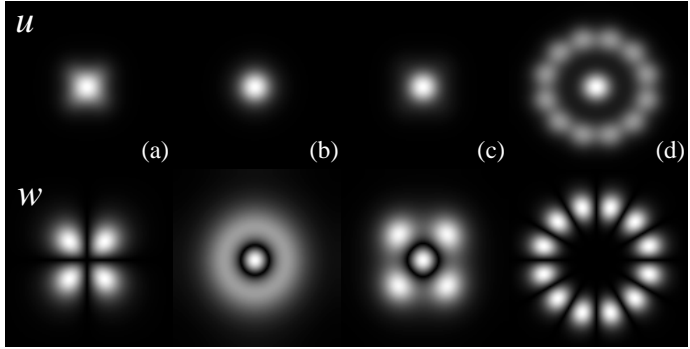,width=0.7\linewidth}
\end{center}
\caption{\label{fig-exotic}
  Different unstable multi--solitonic configurations arising from Eq.
  (\ref{sobolev-nonlineq}) with a change of the initial conditions for fixed
  $\mu_u = -1, \mu_w = -0.3, \kappa=0.5.$}
\end{figure}

Nevertheless the method works equally well for more complicated stationary
solutions.  In Fig. \ref{fig-exotic} we show several exotic solutions that are
obtained by solving Eq.  (\ref{sobolev-nonlineq}) with different initial
conditions. All of those states are dynamically unstable and could have not
been found with traditional minimization methods.

\subsection{Performance and grid refinement}

In this subsection we want to analyze the performance of the method from \S
\ref{excited-search} and to introduce a multigrid--like technique to improve
convergence rates while looking for more accurate solutions.

The idea of the multigrid technique is to use solutions from coarse-grain grids
to calculate better approximations on finer grids \cite{Briggs}. Roughly the
algorithm consists in setting a coarse-grain initial data, solving the equation
or the variational principle with this initial data until the error is small
enough, interpolating the result over a finer grid and iterating using this new
grid until both the error and the spatial discretization of the solution are
the ones we desire.

Since we are already working with Fourier modes over discrete grids the logical
choice for our algorithm is indeed Fourier interpolation. The idea is to use
the expansion from Eq.  (\ref{discrete-fourier}) \emph{outside of the original
  grid}, that is,
\begin{equation}
\psi^{(new)}(x) \equiv \sum_n c_n e^{i\mathbf{k}_n\mathbf{x}},
\quad \forall\mathbf{x}\in\Omega.
\end{equation}
This approximation is then discretized over a finer grid, resulting an
expansion with a larger number of modes
\begin{equation}
\psi^{(new)}_m = \sum_n c^{(new)} e^{i\mathbf{k^{(new)}_n}\mathbf{x}_m},
\end{equation}
which may be used as the starting point for further iteration.

\begin{figure}
\begin{center}
\epsfig{file=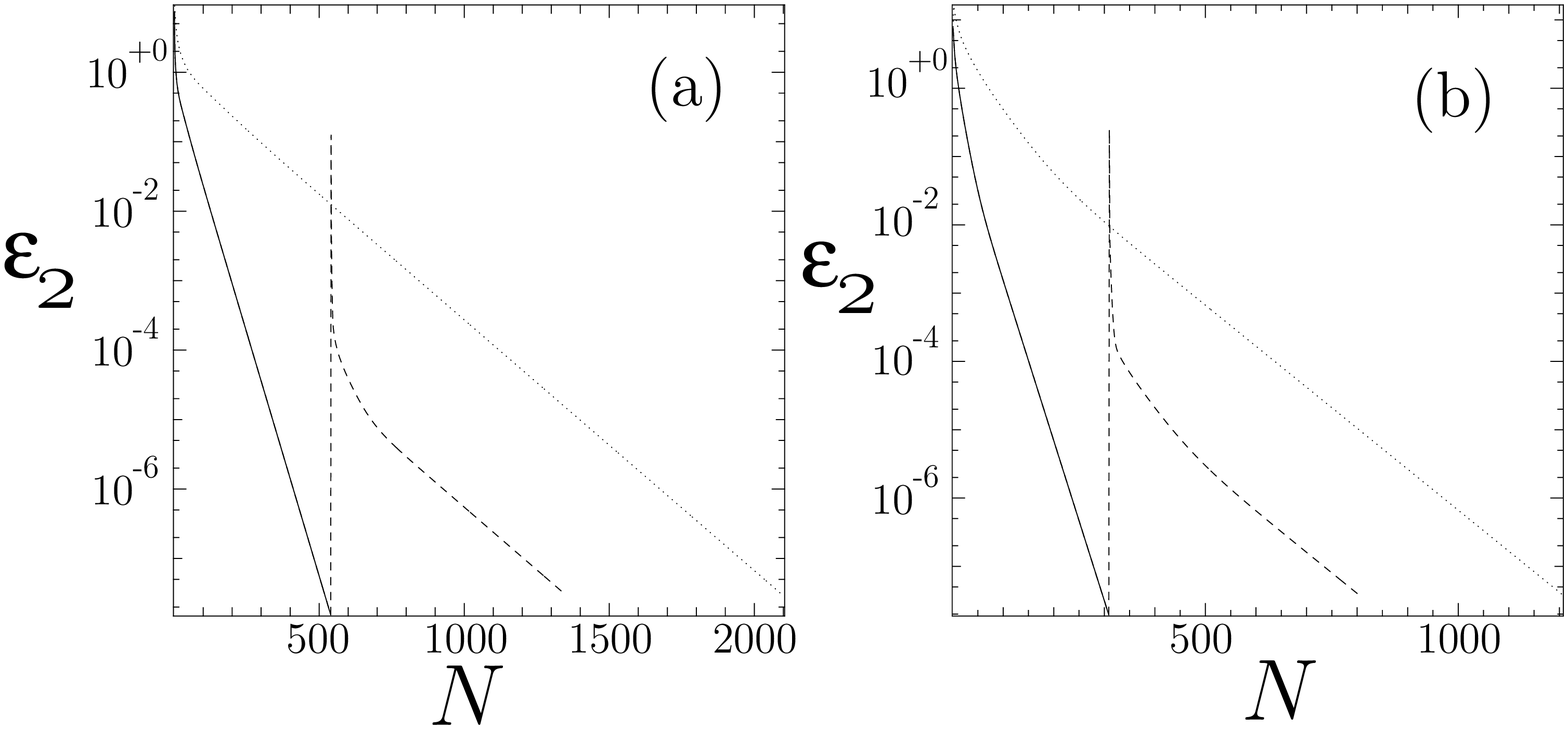,width=10cm}
\end{center}
\caption{\label{fig-perf}
  Dependence of the $L^2$ error on the number of iterations when looking for a (a)
  \emph{vortex-mode vector soliton} or a (b) \emph{dipole-mode vector soliton}.
  We show results for a grid with 32$\times$32 points (solid), 64$\times$64
  points (dotted) and 64$\times$64 points starting from an interpolation of the
  32$\times$32 solution (dashed).}
\end{figure}

In our case we have used only two different grids, one with 32$\times$32 points
(\emph{coarse grid}) and a other with 64$\times$64 points (the \emph{fine
  grid}). We have used two different initial data for our minimization method:
either a pair of Hermite modes which resemble the desired shape, or the
solution of the coarse grid interpolated over the finer grid.  As both the
evolution of the error in Fig.  \ref{fig-perf} and the computation times in
Table \ref{table-perf} show, interpolation saves a significant amount of time.
Indeed, the interpolated solution with only 32$\times$32 Fourier modes is
already a good approximation for the fine grid, as it shows the small change of
the error in Fig. \ref{fig-perf}

\begin{table}
  \vspace{0.3cm}
  \begin{center}
    \begin{tabular}{|c|c|c|c|c|c|c|}   \hline
     & \multicolumn{6}{|c|}{Initial data type }\\ \cline{2-7}
     & \multicolumn{2}{|c|}{Hermite modes} &  Interpolated  &  \multicolumn{2}{|c|}{Hermite modes} &  Interpolated \\
 & \multicolumn{2}{|c|}{for vortex solitons.} & solution  & \multicolumn{2}{|c|}{for dipole solitons.} & solution \\
& & & for vortex & & & for dipole \\ \hline
Grid & 32$\times$32 & 64$\times$64 &   64$\times$64 & 32$\times$32 & 64$\times$64 &   64$\times$64 \\ \hline
  Iterations &  540 & 2101 & 812 &  496 & 1206 & 496 \\ \hline
  Time (s)   &  106 & 2006 & 757 &  60 & 1141 & 466 \\ \hline 
    
   \end{tabular}
  \end{center}
  \vspace{0.3cm}
  \caption{
    \label{table-perf}
    Results for the search of stationary solutions of
    Eq. (\ref{beams-equations1})-(\ref{beams-equations2}), with and without grid refinement.
    Shown are the results for different initial data. The initial data named
    ``interpolated solution'' corresponds to taking as initial data the 
    approximated solution  on the $32\times 32$ grid.}
\end{table}

\section{Conclusions}
\label{sect-final}

In this paper we have shown a numerically efficient way to improve the
convergence of several minimization methods using the so called Sobolev
gradients and applied it to different problems which involve Nonlinear
Schr\"odinger Equations.  We have also proven that these gradients represent a
preconditioning over the traditional definition of gradients on $L^2$. In
Appendix B we suggest a generalization of this method to different vector
spaces and scalar products.

We have presented two different methods for solving our minimization problems:
a traditional one, the \emph{imaginary time evolution}, and a new one,
\emph{the minimization of a nonlinear free energy}.  Both methods have been
shown to be suitable for a Sobolev preconditioning, giving us two new methods
that we call \emph{preconditioned imaginary time evolution} and a
\emph{preconditioned free energy}.

We have implemented all four methods using a discrete Fourier basis and a Fast
Fourier Transform. With these tools we have shown that the Sobolev
preconditioning becomes an inexpensive additional step over existing methods.
The four resulting solvers have been applied to several realistic problems and
in all tests the nonlinear free energy with the Sobolev preconditioning showed
the best convergence rates. Indeed the Sobolev preconditioning has an important
effect on convergence which can be as good as gaining two orders of magnitude
over the traditional techniques.  Furthermore, opposite to what happens with
finite differences \cite{Neu98} the preconditioning may be applied without
significant computational cost.

We have derived a new method for finding excited states of coupled nonlinear
Schr\"odinger equations. This method introduces a new variational principle
which is not based on an energy functional but on finding the zeros of a
nonlinear operator which corresponds to the equation to be solved.  We have
also shown how to improve convergence using Fourier interpolation and a
two-grid method as a source for better initial approximations of the iterative
method. This approach may be extended to more sophisticated multigrid methods.

\appendix

\section{Existence of minimizers for the nonlinear free energy functional.}

\label{free-energy-properties}
In this appendix we want to prove the existence of minimizers for the free
energy (\ref{new-free-energy}). Let us write the free energy functional in the
following form
\begin{equation}
  \label{free-energy-split}
  F(\psi)=\int \bar{\psi}A_\Omega\psi d^{n}x+\int \frac{g}{2}|\psi|^4d^{n}x+\frac{1}{2}\left( N-\lambda \right)^2,
\end{equation}
where $N=\int |\psi|^2d^{n}r$ and
\begin{equation}
  A_\Omega=-\frac{1}{2}\triangle +V(\mathbf{x})-\Omega L_z,\quad \Omega\in[0,1)
\end{equation}
is a positive hermitian operator $A_\Omega\geq \mu _{min}>0.$

Let us also define the following spaces of functions defined over
$\mathbb{R}^n$. We will work in $L^p$ spaces
\begin{equation}
L^p  =\left\{\psi:\, \mathbb{R}^2\rightarrow\mathbb{C} \;/\; \Vert\psi\Vert_p := \left(\int|\psi(\mathbf{x})|^pd^nx\right)^{1/p} < +\infty\right\}.
\end{equation}
Also of interest will be the space of functions over which $A_\Omega$ is well
defined
\begin{equation}
H_\Omega = \left\{\psi:\, \mathbb{R}^2\rightarrow\mathbb{C} \;/\; \Vert\psi\Vert_\Omega := \left(\int\bar{\psi}(\mathbf{x})A_\Omega\psi(\mathbf{x})\right)^{1/2} < +\infty\right\}.
\end{equation}
Let us remark that the dimensionality of the space, $\mathbb{R}^2$ is
important, since for $n\leq 2$ it is easy to show that $H_{\Omega} \subset L^4$.

With the preceding notation we will state the following relevant theorem:

\begin{theorem}
  The free energy functional $F$ given by Eq. (\ref{free-energy-split}) has at
  least one absolute minimum in the set given by the inequalities
  \[
  0 < \Vert\psi\Vert_\Omega^2 \leq \lambda \Vert\psi\Vert_2^2 \leq \lambda(\lambda-\mu_{min}).
  \]
\end{theorem}

\begin{proof}
The first step of the proof will be to show that the domain of $F$
is indeed the whole space $H_\Omega$. Using the positivity of the $A_\Omega$
operator we show that
\begin{equation}
\label{bound-N}
\Vert\psi\Vert_2 \leq \sqrt{\mu_{min}}\Vert\psi\Vert_\Omega,
\end{equation}
which means that $H_\Omega \subset L^2$. We need Sobolev's inequality \cite{Hormander}
\begin{equation}
\Vert\psi\Vert_k \leq \Vert\psi\Vert_2^{1-d} \Vert\nabla\psi\Vert_2^d
\end{equation}
where $d=n/2-n/k$, $n$ is the dimensionality of the space and for us $d=1/2$.
Applying this inequality we obtain the following bound
\begin{equation}
\label{bound-U}
\Vert\psi\Vert_4 \leq \sqrt{\Vert\psi\Vert_2 \Vert\nabla\psi\Vert_2}
\leq \mu_{min}^\frac{1}{4}\Vert\psi\Vert_\Omega,
\end{equation}
which means that $H_\Omega \subset L^4$. Since $\Vert\cdot\Vert_2,\;
\Vert\cdot\Vert_4,\; \Vert\cdot\Vert_\Omega < +\infty$ inside $H_\Omega$ we
conclude that $F$ as given by Eq. (\ref{free-energy-split}) is well defined
over the whole space.

$F$ is also continuous. To prove it let us rewrite the free
energy functional as
\begin{equation}
\label{free-energy-split2}
F(\psi) = \Vert\psi\Vert_\Omega^2 + \frac{g}{2} \Vert\psi\Vert_4^4 + \frac{1}{2}(\Vert\psi\Vert_2^2-\mu)^2,
\end{equation}
and using the bounds (\ref{bound-N}) and (\ref{bound-U}) it is easy to show
that
\begin{equation}
|F(\psi)-F(\xi)| \leq \delta(\varepsilon),\quad \forall \xi\:\Vert\psi-\xi\Vert_\Omega \leq \epsilon
\end{equation}
where the constant $\delta(\varepsilon)$ is given by $\varepsilon$ and
$\Vert\psi\Vert_4$.

We will also need to show that $F$ is coercive
\begin{equation}
\label{coercitivity}
\lim_{\Vert \psi \Vert \rightarrow \infty} \frac{F(\psi)}{\Vert \psi \Vert} \geq \alpha > 0.
\end{equation}
Using Eq. (\ref{free-energy-split2}) one may show that indeed
\begin{equation}
\frac{F(\psi)}{\Vert\psi\Vert_\Omega} \geq \Vert\psi\Vert_\Omega,
\end{equation}
and thus the quotient tends to infinity as the norm grows. The continuity and
coercitivity of $F$ allow to use theorem 1 (section 1.2) of \cite{Direct},
which states the existence of at least one minimizer $\{\inf F(u)\; : \; u\in
X\}$ of $F$ provided it is a weakly lower semicontinuous and coercive
application $F:\,X\rightarrow\mathbb{R}$.

So we know that there is at least one minimum, but we do not know where to look
for it. Let us now show how the $\lambda$ parameter allows us to select
different targets for the minimization problem. To do so we define a real
one-dimensional function for any given direction $\psi \in H_\Omega$
\begin{equation}
  f(k)\equiv F(k\psi).
\end{equation}
This function is nothing but a polynomial over $k$
\begin{equation}
  f(k)=kNa_{\psi}+k^2N^2\frac{u_\psi}{2}+\frac{1}{2}\left( kN-\lambda \right)^2,
\end{equation}
where $a_{\psi}=\int \bar{\psi}A\psi / N$ and $u_{\psi}=\frac{g}{2}\int |\psi |^4/N^2$
are constants that depend only on the precise direction $\psi$.

By differentiating the polynomial and imposing $k=0$ we find that $f'(0) < 0$
at the origin for all possible directions. This means that, as we mentioned
above, our Lagrange ``multiplier'' $\lambda$ allows us to avoid the useless
solution, $\psi = 0$.

Furthermore, we can restrict the location of the minimizer to a surface of
a certain norm. By differentiating $f(k)$ we reach
\begin{equation}
  \label{minimal-norm}
  (a_\psi - \lambda)+ Nk(u_\psi+1) = 0.
\end{equation}
This equation has a single solution which gives us the norm of minimal energy
along the $\psi$-direction
\begin{equation}
  \label{ineq-1}
  N_{min}(\psi) = \max\left\{0,\frac{\lambda-a_\psi}{u_\psi+1}\right\} \leq \lambda - \mu_{min},
\end{equation}
a value which is bounded above by our Lagrange multiplier $\lambda$.

From Eq. (\ref{minimal-norm}) it also follows that for any absolute minimum of
the functional, $\psi_{min}$, the expected value of the $A_\Omega$ operator
must be bounded by the $L^2$ norm
\begin{equation}
  \label{ineq-2}
  \int \bar{\psi}_{min}A\psi_{min} = \Vert\psi_{min}\Vert_\Omega^2 \leq \lambda N(\psi_{min}).
\end{equation}
Otherwise the trivial solution $\psi=0$ would have less energy than $\psi_{min}$.
We can thus, instead of working with an unknown surface, delimit the location
of the minimum to a set which is given by two inequalities, (\ref{ineq-1}) and
(\ref{ineq-2})
\begin{equation}
  W = \{\psi\in H_\Omega\;:\;0 < N(\psi) \leq \lambda-\mu_{min}, \Vert\psi\Vert_\Omega^2 \leq \lambda N(\psi)\}.
\end{equation}
\end{proof}

There are several practical consequences of this theorem. First, it states that
the problem of finding the minima of $E(\psi)$ subject to fixed norm has one
solution, i.e., there exists at least one ground state. Second, but equally
important, it proves our Lagrange penalizer $\frac{1}{2}(N-\mu)^2$ to be
specially well suited for this problem since it avoids both the useless
solution $\psi=0$ and those with too large norm. And finally it gives us a
bounded in which the minimizer must be. Indeed we can extend this result by
proving that we can restrict our search space to a compact superset of $W$.

\begin{theorem}
Any absolute minimum, $\psi$, of the functional $F$ given by Eq. (\ref{free-energy-split})
lays inside the compact set of $L^2$
\[
\bar{U} = \{\psi\in L^2\;:\; \Vert\psi\Vert_\Omega^2 \leq \lambda\Vert\psi\Vert_2^2 \leq \lambda(\lambda-1)\}.
\]
\end{theorem}

\begin{proof}
For the type of spaces that we work with, a set is compact iff it is closed and
we can build an $\varepsilon$-net for any positive number $\varepsilon$.  The
$\varepsilon$-net is a finite set $K_\varepsilon = \{\nu_1,\ldots,\nu_k\}$ such
that for each $\psi \in \bar{U}$ there is an element $\nu \in K_\epsilon$
verifying $\Vert\psi-\nu\Vert_\Omega < \varepsilon$. Thus compactness is
equivalent to the possibility of building an arbitrarily good approximation of
our minimizer using a finite but sufficiently large basis of functions.

It is evident that the set $\bar{U}$ is closed in the subspace $H_\Omega$ of
$L^2$. Let $\{u_n\} \in \bar{U}$ be a convergent sequence and let $u$ be their
limit. Since for each element of the sequence
\begin{equation}
  \Vert u_n \Vert_\Omega \leq \sqrt{\lambda} \Vert u_n \Vert_2 < \lambda(\lambda-1),
\end{equation}
it is also obvious that
\begin{equation}
  \Vert u \Vert_\Omega \leq \sqrt{\lambda} \Vert u \Vert_2 \leq \lambda(\lambda-1),
\end{equation}
also thus the limit belongs to $\bar{U}$.

The compactness of the closed set $\bar{U}$ essentially follows from the fact
that the eigenstates of $A_\Omega$ form a complete basis of $H_\Omega$, and
that the eigenvalues of $A_\Omega$ form a monotonously growing unbounded set of
positive numbers. These eigenstates are of the form
\begin{equation}
  \phi_{n,l} = P_n^l(|\mathbf{x}|)\frac{x_1+ix_2}{|\mathbf{x}|^2}e^{-|\mathbf{x}|^2/2},
\end{equation}
where $P_n^l$ are Laguerre's generalized polynomials. And the corresponding
eigenvalues are
\begin{equation}
  A_\Omega \phi_{n,l} = \mu_{n,l} \phi_{n,l} = (2n+(1-\Omega)l+1)\phi_{n,l},\quad n,l\in\mathbb{N}\cup\{0\}.
\end{equation}

Let us choose any natural number $k$ such that $k+1 > \lambda$. We can split
the whole space as a direct sum $H_\Omega = H_0^{k+1} \oplus H_{k+1}^\infty$
where
\begin{equation}
  H_j^k = \mathrm{lin} \{\phi_{n,l}\;:\;j \leq 2n+l < k\}.
\end{equation}
The important point is that since $H_0^{k+1}$ is isomorph to $\mathbb{R}^m$ for
some natural number $m$, and $\bar{U}_k = \bar{U} \cap H_0^{k+1}$ lays inside
a compact $m$-dimensional ball of radius $\sqrt{\lambda-1}$, then we can find
any $\varepsilon$-net for $\bar{U}_k$. Furthermore, by separating
\begin{equation}
  \psi = \psi_a + \psi_b,\quad \psi_a \in \bar{U}_k, \psi_b \in H_{k+1}^\infty,
\end{equation}
and using the definition of $\bar{U}_k$ we show that the projection of $\psi$
outside of $\bar{U}_k$ can be made arbitrarily small
\begin{equation}
  \Vert \psi_b \Vert_2^2 \leq \frac{\lambda}{k+1} N.
\end{equation}
A direct consequence of this is that for $k > \lambda N / \varepsilon$, a
$\varepsilon$-net of $\bar{U}_k$ is also an $\varepsilon$-net of $\bar{U}$,
which proves the compactness.

Finally, by inspecting the eigenvalues of $A_\Omega$ and using Eq.
(\ref{ineq-1}-\ref{ineq-2}) it is not difficult to see that the absolute minimum of
$F$ must lay in $\bar{U}$.

\end{proof}

\section{Extending the Sobolev gradients}

\label{generalized-gradients} We have shown that redefining the gradient turns
out to be a kind of preconditioning over the original choice of the direction
of descent. Let us assume that our functional has the following form
\begin{equation}
  E(\psi)=\int \bar{\psi}A\psi+f(|\psi|^2,\mathbf{x}),
\end{equation}
where $A$ is a non--negative hermitian operator that may involve some derivatives.
Let us also assume that $H$ is a suitable space equipped with the following
scalar product
\begin{equation}
  \langle \psi,\phi \rangle =\int \bar{\psi}\left( 1+A\right) \phi ,
\end{equation}
which is indeed a scalar product because $\left. A\right|_H\geq 0.$

Let us rewrite the energy functional in the following way
\begin{equation}
  E(\psi)=\langle \psi,\psi\rangle +\int f(|\psi|^2,\mathbf{x})-|\psi|^2.
\end{equation}
The Sobolev gradient in $H$ is
\begin{equation}
  \nabla_AE=\psi+(1+A)^{-1}\left[ \partial_1f-\psi\right] ,
\end{equation}
while the so called ordinary gradient is $\nabla E=A\psi+\partial_1f.$
Hence the preconditioning nature of the method is recovered
\begin{equation}
  \nabla_AE=\left(1+A\right)^{-1}\nabla E.
\end{equation}

We can thus think that the new choice of the scalar product aims at making 
the linear part of the energy functional close to some quadratic form, 
$\langle \psi,B\psi \rangle$, 
such that the new operator $B$ is almost the unity. In our example, indeed,
$B\psi=\psi.$ We that this preconditioning will both enhance the
directions of decreasing energy and also have a smoothing effect on the
nonlinear part. However, the problem is complicated and we know of no
proof that these arguments be of general applicability.

\end{document}